\documentclass[12pt]{article}
\usepackage{amsmath}
\newcommand{\bee}{\begin{eqnarray*}}
\newcommand{\ene}{\end{eqnarray*}}
\newcommand{\beeq}{\begin{equation}}
\newcommand{\eneq}{\end{equation}}
\newtheorem{lem}{Lemma}[section]
\newcommand{\bel}{\begin{lem}}
\newcommand{\enl}{\end{lem}}
\newtheorem{exap}{Example}[section]
\newcommand{\beex}{\begin{exap}}
\newcommand{\enex}{\end{exap}}
\newtheorem{theo}{Theorem}[section]
\newcommand{\beth}{\begin{theo}}
\newcommand{\enth}{\end{theo}}
\newtheorem{prop}{Proposition}[section]
\newcommand{\bep}{\begin{prop}}
\newcommand{\enp}{\end{prop}}
\newtheorem{cor}{Corollary}[section]
\newcommand{\bec}{\begin{cor}}
\newcommand{\enc}{\end{cor}}
\newtheorem{rem}{Remark}[section]
\newcommand{\ber}{\begin{rem}}
\newcommand{\enr}{\end{rem}}
\newtheorem{defi}{Definition}[section]
\newcommand{\bef}{\begin{defi}}
\newcommand{\enf}{\end{defi}}

\usepackage{amsmath}
\usepackage{graphicx}

\begin{document}
\begin{center}
{\large {\bf MLE'S BIAS PATHOLOGY MOTIVATES MCMLE
} }\\
 Yannis G. Yatracos\\
 
 {Cyprus University of Technology, \  30/12/2013}
\end{center}

\date{ }


\begin{center} \vspace{0.25in} {\large Summary} \end{center}


\parbox{4.8in}

{Maximum likelihood estimates  are often biased.
It is shown that this pathology is inherent to the traditional ML estimation  method  for two or more  parameters,
thus motivating  from a different angle the use of MCMLE.
 }\\

\date{ } 

\section{MLE and Bias}

\quad  Various methods have been proposed to reduce the $O(n^{-1})$ term of the asymptotic bias of maximum likelihood estimate (MLE). Firth (1993) observed that 
most  methods are ``corrective'' in character rather than ``preventive'', i.e. the MLE   is first calculated and then 
corrected,  and proposed a preventive approach with  systematic correction of the score equations.
 Bias reduction of MLE's continues to be  a topic of interest
as the current literature indicates; see, for example,  Giles (2012),  Zhang (2013) and the references therein.

In this work it  is shown that bias is inherent to the traditional  ML estimation  method when two or more parameters of a semi-regular  model  
 are estimated. This result
 is confirmed in several examples and  motivates the  use of the preventive Model Corrected MLE
 (Yatracos, 2013),
thus achieving  in these same examples either  partial or 
total bias reduction.
For 
 the Pareto distribution  in particular with both parameters unknown,   the MCMLE $\hat \psi_{MC}$ of the scale parameter $\psi$ improves not only the  bias but also the variance of the 
MLE $\hat \psi.$

\section {The Result-Examples}
 
\quad Let $X$ be a  random vector  from a  parametric model with density $f(x|\theta, \psi),$
 unknown  parameters $ \theta \in R, \ \psi  \in R$ 
and $f$  semi-regular, i.e. at
least the $\psi$-score
$$U_{\psi}(x,\theta, \psi)=\frac{\partial{\log f(x|\theta, \psi)}}{\partial {\psi}}$$
is well defined  (and used to obtain MLE  $\hat \psi$) and  
\begin{equation}
\label{eq:modelreg}
 EU_{\psi}(x, \theta, \psi)=0.
\end{equation}

In the next proposition it is seen  that (\ref{eq:modelreg}) may cause bias of the  MLE $\hat \psi$  because it 
implies often that
$EU_{\psi}(x, \hat  \theta, \psi) \neq 0; \hat \theta$ is the MLE of $\theta.$ Using instead  the score
for the data $Y$ in $U_{\psi}(x, \hat  \theta, \psi)$ this
drawback
is avoided for some models thus motivating from a different  angle the use of MCMLE.

\bep (Bias pathology of  MLE)
\label{p:nsc}
 Let $X$ be a   random vector from the semi-regular parametric model  $f(x|\theta, \psi)$  with $\theta, \ \psi$ both unknown 
and  with score
$ U_{\psi}$
satisfying  (\ref{eq:modelreg}). 
Obtain  MLE  $\hat \theta$ 
either by direct maximization of the likelihood of $X$ or by solving, if it exists, the $\theta$-score  equation 
$$U_{\theta}(x, \theta, \psi)=0.$$
\quad  a) If  
 $ \frac { \partial {U_{\psi}(x, \hat \theta,\psi)} }{ \partial {\psi}}=C$ 
  is  fixed constant,
$ C \neq 0,$
 then $\hat \psi $ is biased estimate of $\psi$ if and only if
\begin{equation}
\label{eq:nsc}
 EU_{\psi}(x, \hat \theta, \psi) \neq  0.
\end{equation} 
Since  (\ref{eq:modelreg}) holds $\hat \psi$ is expected to be biased.\\

 b) If  $ \frac { \partial {U_{\psi}(x, \hat \theta,\psi)} }{ \partial {\psi}}=C(x,\hat \theta, \psi),$ 
$\hat \psi$ is expected more often  to be biased.
\enp

{\bf Proof:} $ a)$  Obtain $\hat \psi$ by  solving   the score equation 
$$U_{\psi}(x,\hat \theta,\psi)=0.$$

Make a  Taylor expansion of $U_{\psi}(x, \hat \theta,\hat  \psi)$ around $ \psi,$ 
\begin{equation}
\label{eq:exp}
U_{\psi}(x, \hat \theta, \hat  \psi)=U_{\psi}(x, \hat \theta,  \psi)+(\hat \psi- \psi)C.
\end{equation}

It follows that 
$$E(\hat \psi- \psi)=-C^{-1} EU_{\psi}(x, \hat  \theta, \psi)  \neq 0$$
 if and only if $EU_{\psi}(x, \hat  \theta, \psi)\neq 0. $ 

$ b)$  Equation (\ref{eq:exp}) remains valid with $C=C(x,\hat \theta, \psi)$ evaluated at $\psi=\psi^*$ between $\psi$ and $\hat \psi.$ Then $\hat \psi$ is biased if and only if
\begin{equation}
\label{eq:expran}
EU_{\psi}(x, \hat \theta, \psi) C^{-1}(x,\hat \theta, \psi^*) \neq 0.
\end{equation}
Make   a second order Taylor approximation of the left  side in (\ref{eq:expran}) around $EU_{\psi}=EU_{\psi}(x, \hat \theta, \psi),\
EC=EC(x,\hat \theta, \psi^*),$
\begin{equation}
\label{eq:Texp}
E\frac{U_{\psi}}{C} \approx \frac{EU_{\psi}}{EC} -\frac{Cov(U_{\psi}, C)}{E^2C}+\frac{Var(C)EU_{\psi}}{E^3C}.
\end{equation}
Whether or not $EU_{\psi}=0,$ (\ref{eq:Texp}) is not expected to vanish.\\

 Proposition \ref{p:nsc} $a)$  shows MLE's  inherent bias pathology since for an ``ideal''   regular model
 (\ref{eq:nsc}) is expected to hold  and thus $\hat \psi$ is biased.
 Proposition \ref{p:nsc} $a)$  holds in all the  examples that follow and MCMLE $\hat \psi_{MC}$ reduces $\hat \psi$'s bias.



\beex
Let $X_1,\ldots, X_n$ be independent  random variables from a normal distribution with mean $\theta$ and variance $\sigma^2$ both unknown.
To use Proposition \ref{p:nsc} $ a)$  for $\sigma^2$ w.l.o.g.  re-parametrize  taking 
$\psi=\sigma^2.$

The log-likelihood  of the sample is
$$-\frac{n}{2} \log \sigma^2-\frac{\sum_{i=1}^n (X_i-\theta)^2}{2\sigma^2}=-\frac{n}{2} \log \psi -\frac{\sum_{i=1}^n (X_i-\theta)^2}{2 \psi},  $$
 the score  equations excluding the constants   are
$$U_{\theta}(X,\theta, \psi)=
\sum_{i=1}^n (X_i-\theta)=0 \rightarrow \hat \theta =\bar X, $$
$$U_{\psi}(X,\hat \theta, \psi)=
-n\psi +\sum_{i=1}^n (X_i-\bar X)^2  =0 \rightarrow \hat \psi= \frac{1}{n}\sum_{i=1}^n (X_i-\bar X)^2 $$
and $\hat \psi$ is biased estimate of $\psi=\sigma^2$ since 
$$\frac{\partial U_{\psi}(X,\hat \theta, \psi)}{\partial \psi}=-n$$
and $EU_{\psi}(X,\hat \theta, \psi) \neq 0.$
\enex

\beex  Let 
$X_1, \cdots, X_n$
be independent random variables from a shifted exponential distribution with  parameters $\theta$ and $\psi (>0)$
both unknown and density
$$f(x,\theta, \psi)=\psi^{-1}e^{-(x-\theta)/\psi}I_{[\theta, \infty)}(x);$$
$I$ denotes the indicator  function. Let $X_{(i)}$ denote  the $i$-th order statistic, $i=1,\ldots, n.$

MLE  $ \hat \theta =X_{(1)}$ and
the score  equation for $\psi,$ after replacing $\theta$ by $X_{(1)},$  is 
$$U_{\psi}(X, \hat \theta, \psi)= -n \psi + \sum_{i=1}^n(X_{(i)}-X_{(1)})=0 \rightarrow \hat \psi=\frac{1}{n} \sum_{i=1}^n(X_{(i)}-X_{(1)}). $$
Since 
$$\frac{\partial U_{\psi}(X,\hat \theta, \psi)}{\partial \psi}=-n$$
 and $EU_{\psi}(X,\hat \theta, \psi) \neq 0,$
 $\hat \psi$ is biased for  $\psi.$
\enex

In the Pareto family example that follows with  parameters $\psi$ and $\theta$ both unknown,
the model corrected MLE, $\hat \psi_{MC},$ of the shape $\psi$ 
reduces  by 50\% the bias  of the MLE $\hat \psi$  and has also smaller variance. With this parametrization $\hat \psi$ is not unbiased
even when $\theta$ is known.  Using the parametrization $\psi=1/\psi^*,$  MLE  $ \hat \psi^*$  is unbiased for $\psi^*$ when $\theta$
is known but when $\theta$ is unknown MCMLE  $\hat \psi^*_{MC}$ is unbiased. 
\beex
 Let 
$X_1, \cdots, X_n$
be independent random variables from a Pareto distribution with density
$$f(x|\theta, \psi)=\psi \theta^{\psi}x^{-(\psi+1)}I_{[\theta, \infty)}(x), \ \psi>0, \ \theta >0;$$
$I$ denotes the indicator function, $n>3.$

The log-likelihood function of the sample is
\begin{equation}
\label{eq:parloglik}
n \log \psi + n \psi \log \theta -(\psi+1) \sum_{i=1}^n \log X_i  + \sum_{i=1}^n \log I_{[\theta, \infty)}(X_i)
\end{equation}
and the MLE estimate of $\theta$ is the smallest observation, $\hat \theta=X_{(1)}.$
The score
$$U_{\psi}(X, \hat \theta, \psi)
= n- \psi \sum_{i=2}^n \log \frac{X_i}{X_{(1)}}$$
and the MLE
 $$\hat \psi=\frac{n}{\sum_{i=2}^n \log \frac{X_i}{X_{(1)}}}.$$ 
Since $Y=\sum_{i=2}^n \log \frac{X_i}{X_{(1)}}$ has a $\Gamma(n-1, \psi)$ distribution (see, e.g,  Baxter, 1980 and references therein)  it follows that $\hat \psi$ is biased and 
$$ \ E\hat \psi-\psi=\frac{2}{n-2}\psi, \ 
Var(\hat \psi)=\frac{n^2}{(n-2)^2(n-3)}\psi^2.$$

The  corrected score based on the data $Y$ is 
$$(n-1)-\psi Y$$
and the MCMLE is
$$\hat \psi_{MC}=\frac{n-1}{\sum_{i=2}^n \log \frac{X_i}{X_{(1)}}},$$
with
$$ E\hat \psi_{MC}-\psi=\frac{1}{n-2}\psi,
 \   Var(\hat \psi_{MC})=\frac{(n-1)^2}{(n-2)^2(n-3)}\psi^2. $$
Observe that $\hat \psi_{MC}$ improves both the bias and the variance of $\hat \psi.$\\

Consider the re-parametrization
$\psi=\frac{1}{\psi^*}.$ The score
$$U_{\psi^*}(X, \hat \theta, \psi^*)
=- n \psi^*+\sum_{i=2}^n \log \frac{X_i}{X_{(1)}}$$
and the MLE $$\hat \psi^*=\frac{\sum_{i=2}^n \log \frac{X_i}{X_{(1)}}}{n}.$$
Proposition \ref{p:nsc} $a)$  holds for
$U_{\psi^*}(X, \hat \theta, \psi^*)$ with $C=-n$ and $EU_{\psi^*}(X, \hat \theta, \psi^*) \neq 0$
indicating that $\hat \psi^*$ is biased.
Using the model from data $Y=\sum_{i=2}^n \log \frac{X_i}{X_{(1)}}$  the corrected score is
$$-(n-1)\psi^*+Y$$ 
and the MCMLE 
$$\hat \psi^*_{MC}=\frac{\sum_{i=2}^n \log \frac{X_i}{X_{(1)}}}{n-1}$$
 is unbiased for  $\psi^*.$ 


\enex



\begin{thebibliography}{xxx99}

\bibitem{xyz}
Baxter, M. A. (1980) Minimum variance unbiased estimation of the parameters of the Pareto distribution. {\em Metrika}, {\bf 27}, 133-138.

\bibitem{xyz}
Firth, D. (1993). Bias reduction of maximum likelihood estimates. {\em Biometrika}, {\bf 80}, 27-38.

\bibitem{xyz}
Giles, D. E. (2012). Bias reduction for the maximum likelihood estimators of the parameters in the
half-logistic distribution. {\em Communications in Statistics – Theory and Methods}, {\bf 41}, 212-222.

 
\bibitem{xyz}
Yatracos, Y. G. (2013) Fisher’s Specification Problem and Model Corrected Maximum Likelihood Estimates (MCMLE).
Submitted for publication.

\bibitem{xyz}
Zhang, J. (2013) Reducing bias of the maximum-likelihood estimation for the truncated Pareto distribution.
{\em Statistics}, {\bf 47}, 792-799.
 

\end{thebibliography}
\end{document}